\documentclass{amsart}
\usepackage{amsthm}
\usepackage{amssymb}
\usepackage{latexsym}
\usepackage{euscript}
\usepackage{amsmath}
\usepackage{amstext}
\usepackage{amsgen}
\usepackage{amsbsy}
\usepackage{amsopn}
\usepackage{amsfonts}
\usepackage{dsfont}
\usepackage[all]{xy}
\usepackage{hyperref}
\usepackage{graphicx}

\title[Some metric properties of spaces of stability conditions]{Some metric properties of\\ spaces of stability conditions}

\author{Jon Woolf}
\date{August 2011}
\thanks{I would like to thank the Newton Institute, Cambridge, for their generous hospitality in April and May 2011, when this paper was completed. I would also like to thank Tom Bridgeland and Arend Bayer for helpful conversations.}

\newtheorem{theorem}{Theorem}[section]
\newtheorem{proposition}[theorem]{Proposition}
\newtheorem{corollary}[theorem]{Corollary}
\newtheorem{lemma}[theorem]{Lemma}

\theoremstyle{definition}

\newtheorem{remark}[theorem]{Remark}

\newcommand{\defn}[1]{\emph{#1}}

\newcommand{\ie}{i.e.\ }

\newcommand{\N}{\mathbb{N}}
\newcommand{\Z}{\mathbb{Z}}
\newcommand{\Q}{\mathbb{Q}}
\newcommand{\R}{\mathbb{R}}
\newcommand{\C}{\mathbb{C}}
\newcommand{\U}{\mathbb{H}}
\renewcommand{\P}{\mathbb{P}}

\newcommand{\im}{\mathrm{im}\,}
\newcommand{\re}{\mathrm{re}\,}

\newcommand{\cat}[1]{\mathcal{#1}}

\newcommand{\id}{\mathrm{id}}

\newcommand{\mor}[2]{{\mathrm{Hom}}(#1,#2)}

\newcommand{\stab}[1]{\mathrm{Stab}(#1)}

\newcommand{\stabo}[1]{\mathrm{Stab}_0(#1)}

\newcommand{\aut}{\mathrm{Aut}}

\newcommand{\dind}{{\overline{d}}}
\newcommand{\dhyp}{d_\mathrm{hyp}}

\newcommand{\mono}{\rightarrowtail}
\newcommand{\epi}{\twoheadrightarrow}

\newdir{ >}{{}*!/-8pt/@{>}}

\begin{document}

\maketitle

\section{Introduction}

The space of locally-finite stability conditions $\stab{\cat{C}}$ on a triangulated category $\cat{C}$ is a (possibly inifinite-dimensional) complex manifold. There is a metric  which induces the topology. For numerically-finite $\cat{C}$, for instance the coherent bounded derived category of a smooth complex projective variety or the bounded derived category of a finite-dimensional algebra, one can consider the finite-dimensional submanifold of numerical stability conditions. We show that under mild conditions this is complete in the metric and describe limiting stability conditions. 
 
The contents are as follows: \S \ref{preliminaries} contains some background material on stability conditions. Section \ref{completeness} contains the proof that a full component of the space of stability conditions whose central charges factor through some finite rank quotient $A$ of the Grothendieck group is complete. Recall that full means that the space of stability conditions is locally homeomorphic to $\mor{A}{\C}$, and not to some proper subspace. The key result is Proposition \ref{limit filtration} which describes the limiting Harder--Narasimhan filtrations. The assumptions are satisified, for instance, for a full component of the space of numerical stability conditions on a numerically-finite triangulated category. 

In  \S \ref{example} we independently verify that the metric on the space of numerical stability conditions on a smooth complex projective curve of genus $\geq 1$ is complete. We compute this metric as follows. There is a natural action of the universal cover $G$ of $GL_2^+\R$ on any space of stability conditions. When the phases of semistable objects are dense for a stability condition $\sigma$, the orbit through $\sigma$ is free and the restriction $d_G$ of the metric to it is independent of $\sigma$ and can be explicitly described, and seen to be complete.  In the case of curves of genus $g\geq 1$ this density of phases condition is satisfied and the action of $G$ is both free and transitive, so that the space of numerical stability conditions is isometric to $(G,d_G)$. The metric $d_G$ is closely related to the hyperbolic metric on the upper half plane --- the universal cover of the conformal linear maps forms a subgroup of $G$ isomorphic to $\C$ and the quotient can be identified with the upper half-plane in such a way that the quotient metric is half the standard hyperbolic metric  ---  see Proposition \ref{explicit metric lemma}.

Section \ref{tilting and hearts of stability conditions} contains some observations about relationships between the hearts of stability conditions. Corollary \ref{tilting corollary} states that hearts of stability conditions in the same component of $\stab{\cat{C}}$ are related by finite sequences of tilts. In the process of proving Theorem \ref{completeness theorem} we obtain a description of the limiting stability condition $\sigma$ of a convergent sequence $\sigma_n$. As a consequence we obtain Corollary \ref{limiting tilt 1}, which states that if the $\sigma_n$ all have the same heart, $\cat{A}$ say, then the heart of $\sigma$ must be a right tilt of $\cat{A}$.

\section{Stability conditions}
\label{preliminaries}

We fix some notation. Let $\cat{C}$ be an additive category. We write $c\in \cat{C}$ to mean $c$ is an object of $\cat{C}$. We will use the term \defn{subcategory} to mean strict, full subcategory. When $\cat{C}$ is triangulated with shift functor $[1]$ exact triangles will be denoted either by $a\to b\to c \to a[1]$ or by a diagram
$$
\xymatrix{
a \ar[rr] && b \ar[dl] \\
& c \ar@{-->}[ul] & 
}
$$
where the dotted arrow denotes a morphism $c \to a[1]$. We will always assume that $\cat{C}$ is essentially small so that isomorphism classes of objects form a set. Given sets $S_i$ of objects for $i\in I$  let $\langle S_i\ |\ i\in I \rangle$ denote the ext-closed subcategory generated by objects isomorphic to an element  in some $S_i$. We use the same notation when the $S_i$ are subcategories of $\cat{C}$.

Let $\cat{C}$ be a triangulated category and $K(\cat{C})$ be its Grothendieck group. A \defn{stability condition} $\sigma=(\mathcal{Z}_\sigma,\mathcal{P}_\sigma)$ on $\cat{C}$ \cite[Definition 1.1]{MR2373143} consists of an additive group homomorphism $\mathcal{Z}_\sigma : K(\cat{C}) \to \C$ and full additive subcategories $\mathcal{P}_\sigma(\varphi)$ of $\cat{C}$ for each $\varphi\in \R$ satisfying
\begin{enumerate}
\item if $c\in \mathcal{P}_\sigma(\varphi)$ then $\mathcal{Z}_\sigma(c) = m(c)\exp(i\pi \varphi)$ where $m(c) \in \R_{>0}$;
\item $\mathcal{P}_\sigma(\varphi+1) = \mathcal{P}_\sigma(\varphi)[1]$ for each $\varphi\in \R$;
\item if $c\in \mathcal{P}_\sigma(\varphi)$ and $c' \in \mathcal{P}_\sigma(\varphi')$ with $\varphi > \varphi'$ then $\mor{c}{c'}=0$;
\item\label{HN filtration} for each nonzero object $c\in \cat{C}$ there is a Harder--Narasimhan filtration, \ie a finite collection of triangles
$$
\xymatrix{
0=c_0 \ar[rr] && c_1 \ar[r] \ar[dl] & \cdots \ar[r] & c_{n-1} \ar[rr] && c_n=c \ar[dl]\\
& b_1\ar@{-->}[ul] &&&& b_n \ar@{-->}[ul]
}
$$
with $b_j \in \mathcal{P}_\sigma(\varphi_j)$ where $\varphi_1 > \cdots > \varphi_n$.
\end{enumerate}
The homomorphism $\mathcal{Z}_\sigma$ is known as the \defn{central charge} and the objects of $\mathcal{P}_\sigma(\varphi)$ are said to be \defn{$\sigma$-semistable of phase $\varphi$}. The objects $b_j$ are known as the \defn{$\sigma$-semistable factors} of $c$. We define  $\varphi_\sigma^+(c)=\varphi_1$ and $\varphi_\sigma^-(c)=\varphi_n$. The \defn{mass}  of $c$ is defined to be $m_\sigma(c)=\sum_{i=1}^n m(b_i)$. 
\begin{lemma}
\label{first phase bounds lemma}
If $a \to b \to c \to a[1]$ is an exact triangle then
$$
\min\{\varphi_\sigma^-(a),\varphi_\sigma^-(c)\}\leq \varphi_\sigma^-(b)\leq \varphi_\sigma^+(b) \leq \max\{\varphi_\sigma^+(a),\varphi_\sigma^+(c)\}
$$
 for any stability condition $\sigma$.
\end{lemma}
\begin{proof}
Note that $\varphi_\sigma^+(b) = \sup \{ t \ | \ \exists\, b' \in \cat{P}_\sigma(t) \ \textrm{with}\ \mor{b'}{b} \neq 0 \}$ and similarly 
$$
\varphi_\sigma^-(b) = \inf \{ t \ | \ \exists\, b' \in \cat{P}_\sigma(t) \ \textrm{with}\ \mor{b}{b'} \neq 0 \}.
$$
Furthermore, both the extreme values are achieved. Suppose $\varphi_\sigma^+(b) > \varphi_\sigma^+(c)$. There is a semistable $b'$ with phase $\varphi^+_\sigma(b)$ and nonzero morphism $b' \to b$. Since $\mor{b'}{c}=0$ for phase reasons this morphism factors through a nonzero morphism to $a$. Hence $\varphi^+_\sigma(a) \geq  \varphi^+_\sigma(b)$ and $\varphi_\sigma^+(b) \leq \max\{\varphi_\sigma^+(a),\varphi_\sigma^+(c)\}$. The other inequality is proved similarly.
\end{proof}
By a \defn{$\sigma$-filtration} of $0\neq c \in \cat{C}$ we mean a filtration as above whose factors are $\sigma$-semistable. The Harder--Narasimhan filtration is one example; in general there will be others not satisfying the decreasing phase of factors condition (which uniquely characterises the Harder--Narasimhan filtration). One consequence of the previous lemma is that the Harder--Narasimhan filtration has minimal phase range amongst all $\sigma$-filtrations in the sense that if a $\sigma$-filtration of $c$ has factors with phases $\varphi_i$ then 
$$
\min\{\varphi_i\} \leq \varphi_\sigma^-(c) \leq \varphi_\sigma^+(c) \leq \max\{\varphi_i\}.
$$

For any interval $I \subset \R$ let $\cat{P}_\sigma(I)=\langle \cat{P}_\sigma(\varphi) \ | \ \varphi\in I \rangle$ be the full subcategory generated by $\sigma$-semistable objects with phases in $I$. For any $t\in \R$ the full subcategory $\cat{P}_\sigma(t,\infty)$ is a $t$-structure. In particular the heart
$$
\cat{P}_\sigma(t,\infty) \cap \cat{P}_\sigma(t,\infty)^\perp[1] = \cat{P}_\sigma(t,\infty) \cap \cat{P}_\sigma(-\infty,t+1] = \cat{P}_\sigma(t,t+1] 
$$
is an abelian subcategory of $\cat{C}$.  Note that $Z_\sigma(a)\neq 0$ for any $0\neq a \in \cat{P}_\sigma(t,t+1]$ so that the phase
$$
\varphi_\sigma(a) = \frac{1}{\pi}\arg Z_\sigma(a) \in (t,t+1]
$$
is well-defined. The semistable objects in $ \cat{P}_\sigma(t,t+1]$ can be characterised as those $b$ such that $a\mono b$ implies $\varphi_\sigma(a)\leq \varphi_\sigma(b)$, or equivalently those for which  $ b \epi c$ implies $\varphi_\sigma(b) \leq \varphi_\sigma(c)$. See \cite[\S5]{MR2373143} for details. The abelian category $\cat{P}_\sigma(0,1]$ is referred to as {\em the} heart of the stability condition $\sigma$.

When the length $|I| =\sup I - \inf I$ of the interval is  $< 1$ the category $\cat{P}_\sigma(I)$ is quasi-abelian, see \cite[Lemma 4.3]{MR2373143}. The strict short exact sequences arise from exact triangles $a \to b \to c \to a[1]$ in $\cat{C}$ for which $a,b$ and $c$ are in  $\cat{P}_\sigma(I)$. A stability condition $\sigma$ is said to be \defn{locally-finite} if given $t\in \R$ we can find $\epsilon >0$ such that the quasi-abelian category $\mathcal{P}_\sigma(t-\epsilon, t+\epsilon)$ is both artinian and noetherian. (Such a category is called finite length in \cite{MR2373143} but we prefer to avoid this terminology to avoid confusion with the {\it a priori} stronger notion of length category introduced in the next section.) 

The set of locally-finite stability conditions can be topologised so that it is a, possibly infinite-dimensional, complex manifold $\stab{\cat{C}}$ referred to as the \defn{space of stability conditions} on $\cat{C}$. For each component of $\stab{\cat{C}}$ there is a linear subspace $V \subset \mor{K(\cat{C})}{\C}$ such that the projection $(\mathcal{Z}, \mathcal{P}) \mapsto \mathcal{Z}$ is a local homeomorphism from that component to $V$ \cite[Theorem 1.2]{MR2373143}. 

The topology on $\stab{\cat{C}}$ arises from the generalised metric
$$
d(\sigma,\tau) = \sup_{0\neq c \in \cat{C}} \max \left\{ | \varphi_\sigma^-(c) - \varphi_\tau^-(c)| , | \varphi_\sigma^+(c) - \varphi_\tau^+(c)|, \left| \log \frac{m_\sigma(c)}{m_\tau(c)}\right| \right\}
$$
which takes values in $[0,\infty]$, see \cite[Proposition 8.1]{MR2373143}. Note that $\{\tau\ | \ d(\sigma,\tau)<\infty\}$ is both open and closed and so is a union of components. Hence $d(\sigma,\tau)=\infty$ implies that $\sigma$ and $\tau$ are in distinct components of $\stab{\cat{C}}$. For the rest of this paper we will loosely refer to $d$ as a metric.

The group $\aut(\cat{C})$ of automorphisms of $\cat{C}$ acts  on the left of $\stab{\cat{C}}$ by isometries  via
$$
(\mathcal{Z}, \mathcal{P}) \mapsto \left(\mathcal{Z}\circ \alpha^{-1}, \alpha\circ \mathcal{P}\right).
$$
There is also a smooth right action of the universal cover $G$ of $GL_2^+\R$. An element $g\in G$ corresponds to a pair $(T_g,\theta_g)$ where $T_g$ is the projection of $g$ to $GL_2^+\R$ under the covering map and $\theta_g:\R \to \R$ is an increasing map with $\theta_g(t+1)=\theta_g(t)+1$ which induces the same map as $T_g$  on the circle $\R/2\Z = \R^2-\{0\} / \R_{>0}$. In these terms the action is given by
$$
(\mathcal{Z}, \mathcal{P}) \mapsto \left(T_g^{-1} \circ \mathcal{Z}, \mathcal{P}\circ \theta_g\right).
$$
(Here we think of the central charge as valued in $\R^2$.) This action preserves the semistable objects, and also preserves the Harder--Narasimhan filtrations of all objects. The subgroup consisting of pairs for which $T$ is conformal is isomorphic to $\C$ with $\lambda\in \C$ acting via 
$$
(\mathcal{Z}, \mathcal{P}) \mapsto \left(\exp(-i\pi\lambda)\mathcal{Z}, \mathcal{P}( \varphi + \re\lambda)\right)
$$
\ie by rotating the phases and rescaling the masses of semistable objects. This action is free and preserves the metric. The action of $1\in \C$ corresponds to the action of the shift automorphism $[1]$. 

The restriction of the metric to an orbit of this $\C$ action is given by $$d(\sigma, \sigma\lambda) = \max\{ |\re \lambda|, \pi|\im \lambda| \}.$$ No sphere in the metric intersects the orbit in a smooth submanifold. Therefore $d$ is not the path metric of a Riemannian metric. \label{not riemannian} Straight lines in $\C$ are geodesics since they have length exactly the distance between their endpoints. By considering geodesic triangles in an orbit we see that $\stab{\cat{C}}$ can never have strictly negative curvature, although it may have non-positive curvature.

Another elementary observation is that the orbits of the $\C$ action are closed. Hence there is an induced quotient metric $\dind$ on $\stab{\cat{C}}/\C$ defined by 
$$
\dind(x\C,y\C) = \inf \{d(x,y\lambda) \ |  \ \lambda \in \C \}.
$$
In \S \ref{example} we will see examples in which this quotient metric arises from a Riemannian metric with constant strictly negative curvature. It would be interesting to see whether  the quotient metric has similar good properties in other examples.

\section{Completeness of spaces of stability conditions}
\label{completeness}

Let $\cat{C}$ be a triangulated category and $K(\cat{C}) \epi A$ a finite rank quotient of its Grothendieck group. Consider the subspace of stability conditions on $\cat{C}$ whose central charge factors through $A$. Let $\stabo{\cat{C}}$ be a component of this subspace. Suppose that $\stabo{\cat{C}}$ is full in the sense of \cite[Definition 4.2]{MR2376815}, \ie that the projection $\stabo{\cat{C}} \to \mor{A}{\C}$ is a local homeomorphism. We prove that such a component is complete in the metric. In fact this follows rather quickly from the results of \cite[\S7,8]{MR2373143}, which show that $\stabo{\cat{C}}$ can be covered by closed balls of fixed radius in the metric, each of which is isometric to a complete metric space.\footnote{I would like to thank Arend Bayer for pointing out this quick argument.}
 However, we give a different proof which has the benefit of providing a useful description of the limiting stability condition.

By assumption stability conditions in $\stabo{\cat{C}}$ are locally-finite. We sharpen this slightly, to show that they are \defn{locally-length} in the sense that given $t\in \R$ we can find $\epsilon >0$ such that the quasi-abelian category $\mathcal{P}_\sigma(t-\epsilon, t+\epsilon)$ is a quasi-abelian length category, or \defn{length category} for short. By this we mean a quasi-abelian category in which every object $a\neq 0$ has a length $l(a)\in \N$ such that any filtration
$$
\cdots \mono a_i \mono a_{i+1} \mono \cdots \mono a
$$
by strict subobjects has length $\leq l(a)$ with equailty for some filtration. (See \cite{MR1779315} for more on quasi-abelian categories.) Clearly a length category is both noetherian and artinian, so that a locally-length stability condition is locally-finite. However, the locally-length condition is {\it a priori} stronger in that it requires the existence of an upper bound on the length of any filtration, not merely that any filtration have finite length. 
\begin{remark}
If $\cat{A}$ is a noetherian and artinian abelian category the Jordan--H\"older theorem asserts that any $0\neq a \in \cat{A}$ has a finite composition series with simple factors, and that all such composition series have the same length. So this definition of length category agrees with the usual one in case $\cat{A}$ is abelian, and any noetherian and artinian abelian category is a length category. Furthermore, for an abelian length category the length descends to a homomorphism $l: K(\cat{A})\to \Z$ defined on the Grothendieck group. However, it is not clear that the analogues hold in the quasi-abelian setting.
\end{remark}

\begin{lemma}
\label{nested length}
Suppose $I \subset J$ are nested intervals of length $<1$. Then $\cat{P}_\sigma(I)$ is a quasi-abelian length category whenever $\cat{P}_\tau(J)$ is. Moreover the length of $a$ in $\cat{P}_\sigma(I)$ is bounded above by its length in $\cat{P}_\tau(J)$.
\end{lemma}
\begin{proof}
A strict monomorphism $\imath:a\mono a'$ in $\cat{P}_\sigma(I)$ fits into a triangle
$a\stackrel{\imath}{\longrightarrow} a' \longrightarrow a'' \longrightarrow a[1]$
with $a,a',a'' \in \cat{P}_\sigma(I) $ and therefore also in $\cat{P}_\tau(J)$. Hence $\imath$ is a strict monomorphism in $\cat{P}_\tau(J)$ too. The result follows.
\end{proof}
The next result is a slightly stronger version of \cite[Lemma 4.4]{MR2376815}; the idea of the proof is the same.
\begin{lemma}
\label{length lemma}
Any $\sigma\in \stabo{\cat{C}}$ is locally-length, indeed $\cat{P}_\sigma(I)$ is a quasi-abelian length category whenever $I$ is a closed interval of length $|I|<1$.
\end{lemma}
\begin{proof}
Suppose $\sigma \in \stabo{\cat{C}}$ and that the central charge $Z_\sigma$ is rational, \ie has image in $\Q[i]$. Then the image is a discrete subgroup of $\C$. Assume, without loss of generality, that $I \subset (0,1)$. Fix $0\neq a \in \cat{P}_\sigma(I)$. Then
$$
K = \left\{ z \in \C \ \left| \ \frac{1}{\pi}\arg z \in I \ \textrm{and} \ \im\, z \leq \im\, Z_\sigma(a) \right. \right\}
$$
is compact, and so contains only finitely many points, say $L$, in the image of $Z_\sigma$. Any filtration of $a$ in $\cat{P}_\sigma(I)$ is mapped under $Z_\sigma$ to a sequence in $K$. Since $Z_\sigma(b)\neq 0$ for $0\neq b \in \cat{P}_\sigma(I)$ the images of the subobjects in the filtration are distinct, and hence the filtration cannot have length greater than $L$. 

In the general case, since $\stabo{\cat{C}}$ is full we can approximate $\sigma$ by a stability condition $\tau$ with rational central charge and such that $\cat{P}_\sigma(I) \subset \cat{P}_\tau(J)$ for some closed interval $J$ of length $<1$. The result follows from Lemma \ref{nested length}.
\end{proof}

Suppose that $\sigma_n$ is a Cauchy sequence in this component. For ease of reading we write $\cat{P}_n, Z_n, \varphi_n,m_n$ for $\cat{P}_{\sigma_n}, Z_{\sigma_n}, \varphi_{\sigma_n}$ and $m_{\sigma_n}$. 
\begin{proposition}
\label{bounded length}
Let $0\neq c\in \cat{C}$. The lengths of the Harder--Narasimhan $\sigma_n$-filtrations of $c$ are bounded. 
\end{proposition}
\begin{proof}
Fix $0< \epsilon < 1/8$ and $\theta \in \R$. Choose $M$ such that $d(\sigma_m,\sigma_n)<\epsilon$ whenever $m,n \geq M$. In particular, for any $n \geq M$,
$$
\cat{P}_n(-\infty,\theta+\epsilon) \subset \cat{P}_M(-\infty,\theta+2\epsilon) 
\quad \textrm{and} \quad 
\cat{P}_n(\theta-\epsilon,\infty) \subset \cat{P}_M(\theta-2\epsilon,\infty).
$$
Let $\tau_n^{< t}$ be the truncation functor right adjoint to the inclusion of $\cat{P}_n(-\infty, t)$ in $\cat{C}$, and so on. Then
\begin{eqnarray*}
\tau_n^{< \theta +\epsilon} \tau_n^{> \theta -\epsilon} c &=& \tau_n^{< \theta +\epsilon} \tau_n^{> \theta -\epsilon} \tau_M^{> \theta -2\epsilon} c \\
&=& \tau_n^{> \theta -\epsilon}  \tau_n^{< \theta +\epsilon} \tau_M^{> \theta -2\epsilon} c \\
&=& \tau_n^{> \theta -\epsilon}  \tau_n^{< \theta +\epsilon} \tau_M^{< \theta +2\epsilon}  \tau_M^{> \theta -2\epsilon} c.
\end{eqnarray*} 
Hence the number of $\sigma_n$-semistable factors of $c$ with phase in $(\theta-\epsilon, \theta+\epsilon)$ is the same as the number of $\sigma_n$-semistable factors of $\tau_M^{< \theta +2\epsilon}  \tau_M^{> \theta -2\epsilon} c$ with phase in $(\theta-\epsilon, \theta+\epsilon)$. However, the Harder--Narasimhan $\sigma_n$-filtration of $\tau_M^{< \theta +2\epsilon}  \tau_M^{> \theta -2\epsilon} c$ is contained within
$$
\cat{P}_n(\theta-3\epsilon,\theta+3\epsilon) \subset  \cat{P}_M(\theta-4\epsilon,\theta+4\epsilon)
$$
and the latter is a length category by Lemma \ref{length lemma}. Therefore there is a uniform (in $n \geq M$) bound on the number of $\sigma_n$-semistable factors of $c$ with phase in $(\theta-\epsilon, \theta+\epsilon)$. Note that $M$ depends only on $\epsilon$, and not on $\theta$. Therefore, covering $[\varphi_M^-(c) -\epsilon, \varphi_M^+(c)+\epsilon]$ by finitely many intervals of the form $(\theta-\epsilon, \theta+\epsilon)$ we obtain a uniform bound on the length of the Harder--Narasimham $\sigma_n$-filtration of $c$ whenever $n \geq M$. 
\end{proof}

For $\theta \in \R$ we define the {\em limiting semistable objects} of phase $\theta$ to be 
$$
\cat{P}(\theta) = \langle 0\neq c\in \cat{C}\ | \ \varphi_n^\pm(c) \to \theta \rangle.
$$
The next lemma shows that limiting semistable objects have non-zero limiting central charge. 
\begin{lemma}
\label{limit semistables}
Let $0\neq c\in \cat{C}$. If $\varphi_n^+(c)-\varphi_n^-(c) \to 0$ as $n \to \infty$ then $Z_n(c) \not \to 0$. 
\end{lemma}
\begin{proof}
If $\varphi_n^+(c)-\varphi_n^-(c) \to 0$ then $m_n(c)- |Z_n(c)| \to 0$. Therefore $Z_n(c)\to 0$ implies $m_n(c)\to 0$, contradicting the assumption that $\sigma_n$ is Cauchy.
\end{proof}

\begin{proposition}
\label{limit filtration}
Each $0\neq c \in \cat{C}$ has a filtration by limiting semistable objects:
$$
\xymatrix{
0=c_0 \ar[rr] && c_1 \ar[r] \ar[dl] & \cdots \ar[r] & c_{r-1} \ar[rr] && c_r=c \ar[dl]\\
& b_1\ar@{-->}[ul] &&&& b_r \ar@{-->}[ul]
}
$$
where the factors $b_i \in \cat{P}(\theta_i)$ with $ \theta_1 >\cdots > \theta_r$. 
\end{proposition}
\begin{proof}
By Proposition~\ref{bounded length} there is some $L\in \N$ such that the Harder--Narasimhan $\sigma_n$-filtration of $c$ has length $\leq L$ for all sufficiently large $n$. We prove the result by induction on $L$. The case $L=1$ is easy, for in this case $c$ is $\sigma_n$-semistable for all sufficiently large $n$ and so $\varphi_n^-(c) - \varphi_n^+(c) \to 0$ (is in fact $0$ for sufficiently large $n$) and $c$ is itself a limiting semistable. Now assume that the result is true for any object whose Harder--Narasimhan $\sigma_n$-filtrations are of length $<L$ for all sufficiently large $n$. 

Since $\sigma_n$ is a Cauchy sequence both $\varphi^+_n(c)$ and $\varphi^-_n(c)$ converge. Suppose that $\varphi_n^+(c)-\varphi_n^-(c) \to \alpha$.  If $\alpha=0$ then $c$ is a limiting semistable and the result holds. So we may assume $\alpha>0$. Choose $0< \epsilon< \alpha / 2(L+1)$. Choose $M$ so that $d(\sigma_m,\sigma_n) < \epsilon$ whenever $m,n\geq M$. For any $n\geq M$ we can find successive semistable factors in the Harder--Narasimhan $\sigma_n$-filtration of $c$  whose phases differ by more than 
$$
\frac{\alpha-2\epsilon}{L} > \frac{2(L+1)\epsilon - 2\epsilon}{L} = 2\epsilon.
$$ 
Splitting the Harder--Narasimhan $\sigma_n$-filtration of $c$ between the above-mentioned factors one obtains a  triangle $c'\to c \to c'' \to c'[1]$ with $c',c''\neq 0$ and such that $\varphi^-_n(c')> \varphi^+_n(c'')$ whenever $n\geq M$. Thus the Harder--Narasimhan $\sigma_n$-filtration of $c$  for any $n\geq M$ is obtained by concatenating those of $c'$ and $c''$. It follows that the Harder--Narasimhan $\sigma_n$-filtrations of $c'$ and $c''$ must have length $<L$ for any $n \geq M$. Hence, by the inductive hypothesis, both $c'$ and $c''$ have filtrations by limiting semistable objects of decreasing phase. Moreover, the minimal phase of a factor of $c'$ is at least as large as the maximal phase of a factor of $c''$. The result follows by concatenating these limiting Harder--Narasimhan filtrations of $c'$ and $c''$.
\end{proof}

\begin{theorem}
\label{completeness theorem}
Let $\stabo{\cat{C}}$ be a full component of the space of locally-finite stability conditions on $\cat{C}$ whose central charges factor through a given finite rank quotient of the Grothendieck group. Then $\stabo{\cat{C}}$ is complete in the natural metric.
\end{theorem}
\begin{proof}
Let $\sigma_n$ be a Cauchy sequence. Then the central charges $Z_n$ converge, say to $Z:K(\cat{C})\to \C$. We claim that $\sigma = (Z,\cat{P})$ is a locally-length stability condition with $\sigma_n \to \sigma$.  It is immediate from the definition of limiting semistable objects that $\cat{P}(\theta+1) = \cat{P}(\theta)[1]$. By Lemma~\ref{limit semistables} if $0\neq c\in \cat{P}(\theta)$ then $Z(c)=m\exp(i\pi \theta)$ for some $m>0$. If $c\in \cat{P}(\theta)$ and $c'\in \cat{P}(\theta')$ with $\theta > \theta'$ then, for sufficiently large $n$, we have $\varphi_n^-(c) > \varphi_n^+(c')$ so that
$\mor{c}{c'}=0$. The existence of Harder--Narasimhan filtrations is guaranteed by Proposition~\ref{limit filtration}. Hence $\sigma$ is a stability condition. If the kernel of $K(\cat{C}) \epi A$ is annihilated by each of the $Z_n$ then it is also annihilated by $Z$, \ie the central charge of $\sigma$ factors through $A$. 

Given $\epsilon >0$ we can choose $M\in \N$ so that $d(\sigma_m,\sigma_n)<\epsilon$ for $m,n\geq M$. Then $\cat{P}(t-\epsilon,t+\epsilon) \subset \cat{P}_m(t-2\epsilon,t+2\epsilon)$ whenever $m\geq M$. Hence $\sigma$ is locally-length by Lemma \ref{nested length}, in particular $\sigma$ is locally-finite. By construction $\sigma_n \to \sigma$, so that $\sigma \in \stabo{\cat{C}}$ and the component is complete.
\end{proof}

\begin{remark}
The condition that $\stabo{\cat{C}}$ is a full component of the space of locally-finite stability conditions on $\cat{C}$ whose central charges factor through a given finite rank quotient of the Grothendieck group enters only in the proof that the lengths of the Harder--Narasimhan $\sigma_n$-filtrations of a fixed object are bounded as $n\to \infty$. Hence Theorem~\ref{completeness theorem} shows that any Cauchy sequence in the space of locally-finite stability conditions for which this is true converges.
\end{remark}

\section{Examples}
\label{example}

In general the metric $d$, and the induced quotient metric on $\stab{\cat{C}}/\C$, are hard to compute. In this section, under the assumption that the phases of $\sigma$-semistables are dense in $\R$, we compute the restricted metric on the orbit $\sigma G$, where $G$ is the universal cover of $GL_2^+\R$, and show that it is independent of $\sigma$. In this case the induced metric on $\sigma G/\C \cong \U$ is half the standard hyperbolic metric on the upper half-plane. This allows us to compute the metric on $\stab{X}$ whenever $X$ is a smooth complex projective curve of genus $\geq 1$, and verify directly that the metrics are complete in these cases. 

\begin{proposition}
\label{explicit metric lemma}
Suppose $\sigma \in \stab{\cat{C}}$ is a stability condition for which the phases of semistable objects are dense in $\R$. Then the $G$ orbit through $\sigma$ is free, the restriction of the metric $d$ to it is independent of $\sigma$ and thus gives a well-defined metric $d_G$ on $G$. This metric is invariant under the action of $G$ by left multiplication and is complete. The induced metric on $\sigma G/\C \cong \U$ is half the hyperbolic metric $\dhyp$ on the upper half-plane.
\end{proposition}
\begin{proof}
Let $g\in G$ correspond to the pair $(T_g,\theta_g)$ where $T_g\in GL_2^+\R$ and $\theta_g:\R \to \R$ is increasing with $\theta_g(t+1)=\theta_g(t)+1$ and induces the same map as $T_g$ on the circle $\R/2\Z = \R^2-\{0\} / \R_{>0}$. Since $\sigma$ and $\sigma g$ have the same semistable objects and Harder--Narasimhan filtrations
\begin{eqnarray*}
d(\sigma g, \sigma) &=& 
\sup_\textrm{semistable $c$} \max \left\{ \left|\varphi_{\sigma g}(c) - \varphi_{\sigma}(c)\right|, \left| \log \frac{m_{\sigma g}(c)}{m_\sigma(c)}\right| \right\}\\
&=&
\sup_\textrm{semistable $c$} \max \left\{ \left| \theta_g(\varphi_{\sigma}(c)) - \varphi_{\sigma}(c)\right|, \left| \log \frac{|T_gZ_{\sigma}(c)|}{|Z_\sigma(c)|}\right| \right\}\\
&=&
\sup_\textrm{semistable $c$} \max \left\{ \left| \theta_g(\varphi_{\sigma}(c)) - \varphi_{\sigma}(c)\right|, \left| \log \left( |T_g v| \right) \right| \right\}
\end{eqnarray*}
where $v = Z_{\sigma}(c) / |Z_\sigma(c)|$ is the unit vector in the direction of $Z_\sigma(c) \in \R^2$. Under the assumption that phases of semistables are dense we thus have
\begin{eqnarray*}
d(\sigma g, \sigma) &=& \max \left\{ \sup_t |\theta_g(t)-t| , \sup_{|v|=1} \left|\log \left( |T_gv| \right) \right| \right\} \\
&=& \max \left\{ || \theta_g -\id||, \log ||T_g||, \log ||T_g^{-1}|| \right\}
\end{eqnarray*}
which depends only on $g$, and not on $\sigma$. Set $\Delta(g)=d(\sigma g, \sigma)$.  It is easy to check that $\Delta(g)=0$ if and only if $g=1$, so that 
$$
d_G(g,h) = d(\sigma g, \sigma h) = d((\sigma h) h^{-1}g , \sigma h) = \Delta(h^{-1}g)
$$
is a metric on $G$, and the orbit $\sigma G$ is free. Note that
 $$
 d_G(fg,fh)=\Delta((fh)^{-1}fg)=\Delta(h^{-1}g) = d_G(g,h)
 $$
 so that $d_G$ is invariant under left multiplication. The explicit description shows that $d_G$ is complete.
 
We have identifications $G/\C \cong GL_2^+\R/\C^* \cong SL_2\R/SO(2)\cong \U$ where the last isomorphism is given by
$$
\left(
\begin{array}{cc}
a & b \\
c& d
\end{array}
\right)
\longmapsto \frac{ai+b}{ci+d}.
$$
Under this identification the left multiplication action of $SL_2\R$ corresponds to its action by M\"obius transformations on $\U$. Hence the induced quotient metric on $\sigma G /\C \cong \U$ is invariant under the M\"obius action of $SL_2\R$. To pin down the metric it therefore suffices to compute the distance between $i$ and $ai$ for real $a\geq 1$. This is given by
$\inf_{\lambda\in \C}\Delta(g\lambda)$ where
$$
T_g = \left(
\begin{array}{cc}
\sqrt{a} & 0 \\
0& 1/\sqrt{a}
\end{array}
\right)
\ \textrm{and}\ 
 \theta_g(t) = \frac{1}{\pi}\tan^{-1}\left( \frac{\tan(\pi t)}{a}\right).
$$
By direct computation we find that $\inf_{\lambda\in \C}\Delta(g\lambda) = \frac{1}{2}\log a = \frac{1}{2}\dhyp(i,ai)$.
\end{proof}
If $X$ is a smooth complex projective curve of genus $\geq 1$ then $\stab{X} \cong G$ 
is a single orbit \cite{MR2335991}. Furthermore, the phases of semistables are dense since, for the standard stability conditon, there are semistable sheaves of any rational slope. (Both of these facts are false for the genus $0$ case.) Therefore the spaces of stability conditions are isometric to $(G,d_G)$ for any such $X$, and in particular are complete. Furthermore, $\stab{X}/\C$ is isometric to the upper half-plane with  metric $\frac{1}{2}\dhyp$.

\section{Tilting and hearts of stability conditions}
\label{tilting and hearts of stability conditions}

A \defn{torsion theory} in an abelian category $\cat{A}$ is a pair $(\cat{T},\cat{F})$ of subcategories such that $\cat{F}\subset \cat{T}^\perp$ and every $a\in \cat{A}$ sits in a short exact sequence $0\to t \to a \to f \to 0$ with $t\in \cat{T}$ and $f\in \cat{F}$, see for example \cite[Definition 1.1]{MR2327478}. In fact $\cat{F}=\cat{T}^\perp$ and $\cat{T}={}^\perp\cat{F}$ so that it is only necessary to specify either $\cat{T}$ or $\cat{F}$. 

Let $\cat{A}$ be the heart of a $t$-structure $\cat{D}$. A torsion theory $\cat{T}$ in $\cat{A}$ determines a new $t$-structure $\langle \cat{D}, \cat{T}[-1] \rangle$, see \cite[Proposition 2.1]{Happel:1996uq}: we say this new $t$-structure is obtained from the original by \defn{left tilting at $\cat{T}$} and denote its heart by $L_\cat{T}\cat{A}$. Explicitly $L_\cat{T}\cat{A} = \langle \cat{F} , \cat{T}[-1] \rangle$. A torsion theory also determines another $t$-structure $\langle \cat{D}^\perp, \cat{F}[1] \rangle^\perp$
by a `double dual' construction. This is the shift by $[1]$ of the other: we say it is obtained by \defn{right tilting} at $\cat{T}$ and denote the new heart $\langle \cat{F}[1], \cat{T} \rangle$ by $R_\cat{T}\cat{A}$.  Left and right tilting are inverse to one another: $\cat{F}$ is a torsion theory in $L_\cat{T}\cat{A}$ and right tilting with respect to this we recover the original heart $\cat{A}$. In the opposite direction, given a $t$-structure $\cat{E}$ with $\cat{D}\subset \cat{E} \subset \cat{D}[-1]$ Beligiannis and Reiten \cite[Theorem 3.1]{MR2327478} show that 
$$
\cat{T}=(\cat{E} \cap \cat{D}^\perp)[1] = \langle a\in\cat{A}\ |\ a=H^1e \ \textrm{for some} \ e\in \cat{E} \rangle
$$ 
determines a torsion theory in $\cat{A}$ such that $L_\cat{T}\cat{A}$ is the heart of $\cat{E}$.

\begin{lemma}
\label{metric and tilting}
If $d(\sigma,\tau)<1/2$ then $\cat{A}_\tau$ can be obtained from $\cat{A}_\sigma$ by performing a left and then a right tilt. (Of course one or both of these tilts may be trivial.)
\end{lemma}
\begin{proof}
Recall that $\cat{A}_\sigma =  \cat{P}_\sigma(0,1]$ is the heart of the $t$-structure $\cat{P}_\sigma(0,\infty)$, and analogously for $\cat{A}_\tau$. Since $d(\sigma,\tau)<1/2$ there are inclusions of $t$-structures
$$
\xymatrix{
& \cat{P}_\tau(0,\infty) \ar[d] & \\
 \cat{P}_\sigma(0,\infty) \ar[r]  &  \cat{P}_\tau\left(-\frac{1}{2},\infty\right) \ar[d] \ar[r]  &  \cat{P}_\sigma(-1,\infty)  \\
 &  \cat{P}_\tau(-1,\infty) &
}
$$
So $\cat{P}_\tau(-1/2,\infty)$ determines torsion theories $\cat{T}$ in $\cat{A}_\sigma$ and $\cat{T}'$ in $\cat{A}_\tau$ with
$$
L_\cat{T} \cat{A}_\sigma = \cat{P}_\tau(-1/2,1/2] = L_\cat{T'}\cat{A}_\tau.
$$
Hence $\cat{A}_\tau = R_{\cat{F}'}L_\cat{T} \cat{A}_\sigma$.
\end{proof}
\begin{corollary}
\label{tilting corollary}
If $\sigma$ and $\tau$ are in the same component of $\stab{\cat{C}}$ then the hearts $\cat{A}_\sigma$ and $\cat{A}_\tau$ are related by a finite sequence of left and right tilts.
\end{corollary}
\begin{proof}
The connected components of $\stab{\cat{C}}$ are the path components. Choose a path from $\sigma$ to $\tau$, cover it with finitely many balls of diameter $<1/2$ and apply Lemma~\ref{metric and tilting}.
\end{proof}

For a heart $\cat{A}$ let $U(\cat{A}) = \{\sigma \in \stab{\cat{C}} \ | \ \cat{A}_\sigma=\cat{A} \}$. Note that $U(\cat{A})$ may be empty and need not be either open or closed. The next result describes the hearts of stability conditions in the closure of $U(\cat{A})$.
\begin{corollary}
\label{limiting tilt 1}
Suppose $\sigma_n$ is a sequence in $U(\cat{A})$ with $\sigma_n \to \sigma$. Then $\cat{A}_\sigma = R_\cat{T}\cat{A}$ for the torsion theory with  
$$
\cat{T} = \langle 0\neq a \in \cat{A} \ | \ \varphi_n^-(a)\not \to 0 \rangle
$$ 
and $\cat{F} = \langle 0\neq a \in \cat{A} \ | \ \varphi_n^+(a) \to 0 \rangle$.
\end{corollary}
\begin{proof}
It is easy to check that $\cat{T}$ is a torsion theory; the required short exact sequences arise from the description of the Harder--Narasimhan $\sigma$-filtration as a limiting filtration. The description of the limiting semistable objects shows that $R_\cat{T}\cat{A} \subset \cat{A}_\sigma$. Hence they are equal as hearts of distinct non-degenerate $t$-structures cannot be nested. 
\end{proof}

It is not always the case that a torsion theory has the form $\langle 0\neq a \in \cat{A} \ | \ \varphi_n^-(a) \not \to 0 \rangle$ for some Cauchy sequence $\sigma_n$ of stability conditions, and it is not true that $\cat{B}=R_\cat{T}\cat{A}$ implies $\overline{U(\cat{A})}\cap  U(\cat{B}) \neq \emptyset$. For example, let $\cat{A}$ be the category of representations of the Kronecker quiver. This is well-known to be derived equivalent to the coherent derived category of $\P^1$ and (after fixing an appropriate equivalence) $\cat{A}=L_\cat{T}{\rm Coh}(\P^1)$ where $\cat{T}={}^\perp\langle \mathcal{O}(d) \ |\ d<0 \rangle$. However there is no sequence $\sigma_n$ in $U(\cat{A})$ with $\sigma_n \to \sigma$ where $\cat{A}_\sigma={\rm Coh}(\P^1)$; see \cite[\S 3.2]{Woolf:2010fk} for details.

\end{document}